\documentstyle[11pt,psfig]{article}

\title{ANALYSIS OF THE VIBRATIONAL MODE SPECTRUM OF A LINEAR CHAIN WITH SPATIALLY EXPONENTIAL PROPERTIES}

\author{ {\sl Thomas M.  Michelitsch $^{1}$\footnote{Corresponding author, Email: michel@lmm.jussieu.fr}, G\'{e}rard A. Maugin $^{1}$, Andrzej F. Nowakowski$^2$, Franck C. G. A. Nicolleau$^2$} \\[10pt]
 $^1$ Institut Jean le Rond d'Alembert \\ CNRS UMR 7190 \\
 Universit\'{e} Pierre et Marie Curie, Paris 6\\ FRANCE
 \\[5pt]
 \\
$^2$ Department of Mechanical Engineering\\
University of Sheffield, UK}

\date{ }

\begin{document}

\maketitle

\paragraph{\it International Journal of Engineering Science 47(2), 209-220 (2009). \\
$\rm doi:10.1016/j.ijengsci.2008.08.011$ \\
Dedicated to Valery Levin's 70th Birthday}

\paragraph{ABSTRACT}

We deduce the dynamic frequency-domain-lattice Green's function of a linear chain with  properties (masses and next-neighbor spring constants) of exponential spatial dependence.
We analyze the system as discrete chain as well as the continuous limiting case which represents an elastic 1D exponentially graded material.
The discrete model yields closed form expressions for the $N \times N$ Green's function for an arbitrary number $N=2,..,\infty$ of particles of the chain. Utilizing this Green's function yields an explicit expression for the vibrational mode density. Despite of its simplicity the model reflects some characteristics of the dynamics of a 1D exponentially graded elastic material. As a special case the well-known expressions for the Green's function and oscillator density of the homogeneous linear chain are contained in the model.
The width of the frequency band is determined by the grading parameter which characterizes the exponential spatial dependence of the properties.
In the limiting case of large grading parameter, the frequency band is localized around a single finite frequency where the band width tends to zero inversely with the grading parameter. In the continuum limit the discrete Green's function recovers the Green's function of the continuous equation of motion which takes in the time domain the form of a Klein-Gordon equation.


\section{INTRODUCTION}

There is an increasing interest in many dynamic applications to control the frequency band
of vibrational modes.
Extensive models dealing with how point defects in crystalline lattices affect the density of vibrational modes by using lattice Green's functions \cite{dederichs}.
As a rule in these models the defects are conceived as sufficiently "small" perturbations
of a regular lattice
and treated in the framework of a perturbation calculus where the lattice is considered in the harmonic approximation \cite{maradudin,montroll,mougios}. All these models have in common that they consider the effect of localized defects on the vibrational spectrum resulting in the appearance of "soft modes" which can be verified by their effects on the low temperature
behavior of the specific heat \cite{mougios,zeller}.

However, little is known about the vibrational mode density when the lattice has "strongly" spatially varying properties which cannot be conceived as "small perturbations" with respect to an homogeneous reference lattice.
In recent decades the study of spatially inhomogeneous materials has
become of increasing interest, in particular due to their increasing importance in new applications in so called Smart Materials and as functionally graded materials (FGM).
However their technological exploitation still is limited in view of the principal difficulty to set up analytic continuum models. A principal reason for this is the fact that even for their simplest variants, the FGM with exponential properties,
continuous Green's functions are not available in explicit forms.
This is true for static Green's functions \cite{1} and more than ever for the dynamic framework such as as for problems of heat conduction \cite{2} in two and three dimensions.

The present paper aims at establishing a simple discrete one-dimensional linear chain model for exponentially graded material. The exponential spatial dependence of the properties, i.e. of the masses and spring constants is assumed such that their ratios remain spatially independent. The frequency-domain dynamic lattice Green's function of this 1D system is determined in closed form for the entire range of frequencies, that is within the band of eigenfrequencies and outside the band of frequencies.
The expression deduced for the $N\times N$ Green's function matrix is {\it explicit} for an arbitrary number $N>1$ of particles in the chain. Utilizing this Green's function gives also an exact expression the density of oscillation modes. It turns out that variation of the frequency band as well as the lowest and highest frequency can be controlled by the grading parameter. In spite of the simplicity of the model, the latter provides some of the essential features of one dimensional FGM.

\section{The Discrete Model}
We analyze a linear chain consisting of $N>>1$ mass particles which are
assumed to be harmonically connected to their next neighbors.
Any particle in the chain $p$ ($p=0,..,N-1$) is assumed to have one degree of freedom characterized by the displacement $u_p(t)$, where $t$ indicates the time coordinate.
The Hamiltonian of this system can be written as
\begin{equation}
\label{hamiltonnth}
H=\frac{m_0}{2}\sum_{p=0}^{N-1}\xi^{2p}\left\{\dot{u}_p^2+\omega_0^2(u_p-u_{p+1})^2\right\}
\end{equation}

The grading parameter $\xi$ characterize the exponential spatial dependence of the masses $m_p=m_0\xi^{2p}$ and of the spring constants $f_p=m_p\omega_0^2$.
We assume that $\xi>0$, where $\xi=1$ represents the "trivial" case of a homogeneous linear chain. A schematic representation of the chain is drawn in Fig. 1.

Our goal is to deduce the $N \times N$ frequency domain lattice Green's function containing the complete dynamic information such as the vibrational mode spectrum.
The Hamiltonian equations of motion of the above system (\ref{hamiltonnth}) are obtained from

\begin{equation}
\label{eq0}
m_0\xi^{2p}\,\,\ddot{u}_p=-\frac{\partial H}{\partial u_p}
\end{equation}
where absence of external forces is assumed. This equation reads

\begin{equation}
\label{eq1}
\xi^{2p}\,\,\ddot{u}_p=-\omega_0^2\left\{\xi^{2p}(u_p-u_{p+1})+\xi^{2p-2}(u_p-u_{p-1})\right\} \,,\hspace{1cm} p=0,..,N-1
\end{equation}
where strictly speaking this equation is not defined at the boundaries $p=0,N-1$ without imposing boundary conditions, i.e. assumptions about the fictive displacements $u_{-1}, u_N$. In the case of a homogeneous chain (represented by $\xi=1$) it is convenient to impose periodic boundary conditions. Let us for a moment disregard boundary conditions to re-introduce them later in a quasi "natural" way.
At present (\ref{eq1}) can be rewritten more compactly in matrix form as

\begin{equation}
\label{newt1}
 \ddot{{\bf u}}=-{\bf \Lambda}^{-1}\cdot{\bf K}\cdot{\bf u}
\end{equation}

Here we introduced the (dimensionless) $N\times N$ mass-matrix ${\bf \Lambda}$
\begin{equation}
\label{mass1}
\Lambda_{pq}=\delta_{pq}\xi^{2p}
\end{equation}
where $\delta_{pq}$ denotes the Kronecker symbol (no sum over $p$ is performed) and
the $N\times N$ stiffness matrix ${\bf K}$ which is given by
\begin{equation}
\label{Kstiffmat}
K_{pq}=\xi^{2p}\omega_0^2\left(\delta_{pq}(1+\frac{1}{\xi^2})-\delta_{p+1,\,q}-\frac{1}{\xi^2}\delta_{p-1,\,q}\right)
\end{equation}

Whether or not ${\bf K}$ is positive definite or positive semi-definite, i.e. whether the
eigenfrequency zero exists which refers to a uniform translation $u_p=at+b$ $\forall p$ of the chain, depends on whether or not such uniform motion is compatible with the boundary conditions.
We have to note that the matrix ${\bf \Lambda}^{-1}{\bf K}$ of (\ref{newt1}) is {\it non-symmetric}, i.e. non Hermitean. Hence its eigenmodes are mutually non-orthogonal.

To determine the eigenmodes and eigenvalues of the "dynamic matrix" ${\bf \Lambda}^{-1}{\bf K}$ it is convenient to introduce an auxiliary vector ${\bf y}$ defined by\footnote{By this substitution we transform the problem into a self-adjoint one with respect to the new
variable ${\bf y}$.}

\begin{equation}
\label{subst1}
{\bf y}={\bf \Lambda}^{\frac{1}{2}}\cdot{\bf u}
\end{equation}
having the Cartesian components $y_p=\xi^p\,u_p$ (no sum over $p$). The equation of motion for ${\bf y}$ assumes then the form

\begin{equation}
\label{newton2}
 \ddot{{\bf y}} =-{\bf L}\cdot{\bf y}\,,\hspace{3cm} {\bf L}={\bf \Lambda}^{-\frac{1}{2}}\,{\bf K}\,{\bf \Lambda}^{-\frac{1}{2}}
\end{equation}

Matrix ${\bf L}$ has the components $L_{pq}=\xi^{-(p+q)}K_{pq}$ (no sum over $p,q$) with

\begin{equation}
\label{matrixL1}
L_{pq}=\frac{\omega_0^2}{\xi^2}\left(\delta_{pq}(1+\xi^2)-\xi(\delta_{p+1,\,q}+\delta_{p-1,\,q})\right)
\end{equation}

We observe that $L_{pq}=L_{qp}$ is a {\it symmetric}, Hermitean (self-adjoint) matrix which coincides in the degenerate case $\xi=1$ with $K_{pq}$ representing the homogeneous linear chain. Hence the Hamiltonian (\ref{hamiltonnth}) can be rewritten as

\begin{equation}
\label{hamil}
H=\frac{m_0}{2}\left(\sum_{p=0}^{N-1}\dot{y}_p^2+\sum_{p,q=0}^{N-1}L_{pq}y_py_q\right)
\end{equation}

The rearrangement of terms which lead to expression (\ref{hamil}) is justified when we assume
periodic boundary conditions for $y_p$, namely
\begin{equation}
\label{peroidboundc}
y_p(t)=y_{p+N}(t)
\end{equation}

Let us further take into account that

\begin{equation}
\label{matL2}
{\bf L}={\bf \Lambda}^{\frac{1}{2}}\,\left({\bf \Lambda}^{-1}\,{\bf K}\right)\,{\bf \Lambda}^{-\frac{1}{2}}
\end{equation}

From this relation follows that ${\bf L}$ has the same spectrum of eigenvalues as the dynamic matrix ${\bf \Lambda}^{-1}{\bf K}$ which characterizes the motion of the
"true" displacements $u_p=\xi^{-p}y_p$. Hence we can determine the frequency spectrum of ${\bf \Lambda}^{-1}{\bf K}$ by just determining the spectrum of ${\bf L}$.
Due to the symmetry of ${\bf L}$ it has has mutually orthogonal eigenvectors (unlike ${\bf \Lambda}^{-1}{\bf K}$).

Equation (\ref{matrixL1}) suggests that the eigenvectors ${\bf v}$ of ${\bf L}$ have components of the form

\begin{equation}
\label{eigen1}
v_p=C \mu^p
\end{equation}
with $C$ being a constant and the admissible values of $\mu$ are yet to be determined. Then we have
\begin{equation}
\label{eigveceq}
\sum_{q=0}^{N-1}L_{pq}v_q=C \mu^p \,\frac{\omega_0^2}{\xi^2}\left\{(1+\xi^2)-\xi(\mu+\mu^{-1})\right\}= \lambda v_p
\end{equation}
From this relation we find
that to these eigenvectors correspond obviously the eigenvalues
\begin{equation}
\label{eigenvals1}
\lambda= \frac{\omega_0^2}{\xi^2}\left\{(1+\xi^2)-\xi(\mu+\mu^{-1})\right\}=
\frac{\omega_0^2}{\xi^2}\left\{(\xi-\mu)(\xi-\mu^{-1})\right\}
\end{equation}

Taking into account the periodic boundary conditions (\ref{peroidboundc}) which must be
reflected by the eigenvectors to fulfill these boundary condition, namely $v_p=v_{p+N}$, leads to the condition for $\mu$ as
\begin{equation}
\label{condmu}
\mu^N=1
\end{equation}

Hence the admissible $\mu$ are just the $N$ complex unity roots given by

\begin{equation}
\label{unityroots}
\mu_m=e^{ik_m}\,,  \hspace{2cm}   k_m= \frac{2\pi m}{N}, \,\,\,\,\,m=0,..,N-1
\end{equation}

By imposing periodic boundary conditions (\ref{peroidboundc}) we define at the same time
${\bf L}$ in (\ref{matrixL1}) at the boundaries for $p=0,N-1$ by its periodicity condition

\begin{equation}
\label{Lbound}
L_{pq}=L_{p+N,q}, \,\, p=0,..,N-1
\end{equation}
which holds for both indices $p$ and $q$ because of its symmetry $L_{pq}=L_{qp}$.
With (\ref{peroidboundc}) both ${\bf L}$ and ${\bf K}$ become completely defined $N\times N$ matrices, i.e. for $\xi\neq 1$ {\it positive definite}\footnote{Since the uniform translation $u_p(t)=at+b$ is not compatible with
(\ref{peroidboundc})} and for $\xi=1$ positive semi-definite.
The $N$ normalized (Bloch) eigenvectors ${\bf v}^{(m)}$ ($m=0,..,N-1$) of ${\bf L}$ have the components

\begin{equation}
\label{blocheig}
v_p^{(m)}=\frac{1}{\sqrt{N}} e^{ik_mp}
\end{equation}
fulfilling the condition of ortho-normality
\begin{equation}
\label{blochoth}
\sum_{p=0}^{N-1}v_p^{(m)}v_p^{(n)*}=\delta_{mn}
\end{equation}
and the completeness condition

\begin{equation}
\label{blochothcmp}
\sum_{m=0}^{N-1}v_p^{(m)}v_q^{(m)*}=\delta_{pq}
\end{equation}

With $k_m=\frac{2\pi}{N}m$ the explicit $N$ eigenvalues read

\begin{equation}
\label{eigenvals}
\lambda_m= \omega_m^2=\frac{\omega_0^2}{\xi^2}\left\{(\xi-\cos{k_m})^2+\sin^2{k_m}\right\}, \,\, m=0,..N-1
\end{equation}
with the dispersion relation $\omega_m=\omega(k_m)$. For $\xi\neq 1$ it follows that $0<\frac{\omega_0^2}{\xi^2}(1-\xi)^2\leq \lambda_m\leq \frac{\omega_0^2}{\xi^2}(1+\xi)^2$. The frequency $\omega_0=0$ occurs for $k_0=0$ only in the degenerate case of a homogeneous chain when $\xi=1$. By introducing the periodic boundary condition (\ref{peroidboundc}) matrix ${\bf L}$ has for $\xi \neq 0$
no eigenvalue equal to zero. This is due to the fact that (\ref{peroidboundc}) does not admit
uniform translations in all components $u_p=\xi^{-p}y_p$ which would correspond to an eigenvector equal to zero. The following further observations can be made:
\newline\newline
\noindent (I) Limiting case $\xi >>1$: This limiting case yields for all eigenfrequencies
$\lim_{\xi\rightarrow\infty}\omega_m=\omega_0$ where the band width tends to zero as
$2\frac{\omega_0}{\xi}$. By increasing the grading parameter one can "squeeze" the
eigenfrequencies $\omega_m$ all being localized within the band $\omega_0(1-\frac{1}{\xi})\leq\omega_m\leq\omega_0(1+\frac{1}{\xi})$.
\newline\newline
\noindent (II) Limiting case $\xi\rightarrow 0$:  $\lim_{\xi\rightarrow 0}\omega_m=\frac{\omega_0}{\xi}\rightarrow\infty$ where the frequencies are localized within
the band $\omega_0(\frac{1}{\xi}-1)\leq \omega_m \leq \omega_0(\frac{1}{\xi}+1)$
of constant width $2\omega_0$. By decreasing the grading parameter we can "drive" the eigenfrequencies to arbitrarily high values where the band width remains constant $2\omega_0$.
\newline\newline
The possibility to control the band of eigenfrequencies by the grading parameter may open interesting dynamic applications. From (\ref{matL2}) it follows that ${\bf \Lambda}^{-1}{\bf K}$ has the same eigenvalues
(\ref{eigenvals}) where its eigenvectors have the components $\xi^{-p}v_p^{(m)}=\frac{\xi^{-p}}{\sqrt{N}} e^{ik_mp}$ ($m=0,..,N-1$). We observe that the $p^{th}$ components of {\it each} eigenvector ${\bf v}^m$ are
scaled with the same scaling factor $\xi^{-p}$ which is therefore true for any admissible displacement
$u_p$ which thus takes the modal expansion
\begin{equation}
\label{displ}
u_p(t)= \xi^{-p}y_p(t)= \frac{\xi^{-p}}{\sqrt{N}}\sum_{m=0}^{N-1}\left(A_m e^{i(k_mp-\omega_mt)}+B_m e^{i(k_mp+\omega_mt)}\right)
\end{equation}
Equation (\ref{displ}) represents the general solution of the homogeneous\footnote{"Homogeneous" due to the absence of external forces.} system (\ref{newt1}) with
$2N$ arbitrary constant coefficients $A_m, B_m$  which can be uniquely fixed to satisfy any 2N initial data $u_p(t=0),\dot{u}_p(t=0)$. The motion associated with (\ref{displ}) is hence characterized by
the Bloch waves $y_p(t)$ being spatially exponentially deformed by the factor $\xi^{-p}$.
From (\ref{displ}) we observe the scaling relation $u_{p+N}(t)=
\xi^{-N}u_p(t)$.

For the following deduction of the spectral oscillation mode density density it is convenient to employ the apparatus of dynamic Green's functions. Let us first define the time-domain Green's function ${\bf {\cal G}}(t)$ describing the motion of ${\bf y}(t)$. To this end we consider the {\it inhomogeneous} differential equation system

\begin{equation}
\label{gfydef}
\left\{(\frac{d}{dt}+\epsilon)^2{\bf 1}+{\bf L}\right\}\cdot{\bf y}={\bf b}(t)
\end{equation}
where $\epsilon>0$ denotes an infinitesimal positive damping constant, ${\bf 1}=(\delta_{pq})$
the $N\times N$ unity matrix, and ${\bf b}(t)$ is an arbitrary source term which determines the response ${\bf y}(t)$ uniquely. Expressing the solution ${\bf y}(t)$ in terms of a convolution yields
a representation of the form

\begin{equation}
\label{convol1}
{\bf y}(t)=\int_{-\infty}^{\infty}{\bf {\cal G}}(t-t')\cdot {\bf b}(t'){\rm d}t'
\end{equation}
where the kernel ${\bf {\cal G}}$ is the (time domain) Green's function of (\ref{gfydef}). By taking into account the trivial identity

\begin{equation}
\label{trivid}
{\bf b}(t)= \int_{-\infty}^{\infty} \delta(t-t'){\bf 1}\cdot{\bf b}(t'){\rm d}t'
\end{equation}
thus (\ref{gfydef}) requires that
\begin{equation}
\label{greendf1}
\left\{{\bf 1}(\frac{d}{dt}+\epsilon)^2+{\bf L}\right\}{\bf {\cal G}}(t)= \delta(t){\bf 1}
\end{equation}
which defines the Green's function ${\bf {\cal G}}(t)$. By introducing the Fourier transforms

\begin{equation}
\label{foutragf1}
{\bf {\cal G}}(t)=\frac{1}{2\pi}\int_{-\infty}^{\infty}{\bf {\hat {\cal G}}}(\omega)e^{-i\omega t}{\rm d}\omega
\end{equation}
and

\begin{equation}
\label{foutdelta}
\delta(t)=\frac{1}{2\pi}\int_{-\infty}^{\infty}e^{-i\omega t}{\rm d}\omega
\end{equation}
we can transform (\ref{greendf1}) into the frequency domain to determine
the dynamic Green's function ${\bf {\hat {\cal G}}}(\omega)$ in the frequency domain (which is nothing but the Fourier transform of the time domain Green's function ${\bf {\cal G}}(t)$) to arrive at

\begin{equation}
\label{greendf2}
{\bf {\hat {\cal G}}}(\omega)= \left\{({\bf L}-{\bf 1}(\omega+i\epsilon)^2\right\}^{-1}
\end{equation}
where ${\hat {\cal G}}_{pq}={\hat {\cal G}}_{qp}$ has the same symmetries as ${\bf L}$
and contains the complete dynamic information on the system, especially on the
oscillation spectrum (density of normal oscillators) which we will deduce below.
Note that we introduced an infinitesimal damping constant $\epsilon >0$ which guarantees that the Fourier representation (\ref{foutragf1}) of the time-domain Green's
function ${\bf {\hat {\cal G}}}(\omega)$ remains well defined, i.e. ${\bf L}-{\bf 1}(\omega+i\epsilon)^2$ remains invertible along the entire real $\omega$-axis,
especially also at the eigenfrequencies $\omega=\omega_m$. In the following any expression of the form $f(\omega+i\epsilon)$ is conceived as its limiting case $\epsilon\rightarrow 0+$ where the zero is approached from the positive side $\epsilon>0$. Inverting (\ref{greendf1})
requires such an infinitesimal damping constant which regularizes the problem and
guarantees the time-domain Green's function to be causal \cite{michelwang} (i.e. ${\bf {\cal G}}(t)=\Theta(t)\,{\bf {\cal G}}(t)$ where $\Theta(t)$ denotes the Heaviside unit step function).

Before we evaluate (\ref{greendf2}) explicitly it is worthy to consider its relation to the (non-symmetric) Green's function
of the "true" displacements ${\bf u}(t)$ which is defined by

\begin{equation}
\label{gftrue}
{\hat {\bf G}}(\omega)=\left[{\bf \Lambda}^{-1}\,{\bf K}-{\bf 1}(\omega+i\epsilon)^2\right]^{-1}
\end{equation}

The transformation between ${\hat {\bf G}}$ and ${\bf {\hat {\cal G}}}$ is easily gained by rewriting ({\ref{gftrue}) in the form
\begin{equation}
\label{gftrue2}
{\hat {\bf G}}(\omega)=\left[{\bf \Lambda}^{-\frac{1}{2}}\left({\bf \Lambda}^{-\frac{1}{2}}\,{\bf K}{\bf \Lambda}^{-\frac{1}{2}}-{\bf 1}(\omega+i\epsilon)^2\right){\bf \Lambda}^{\frac{1}{2}}\right]^{-1}
\end{equation}
or
\begin{equation}
\label{gftrue3}
{\hat {\bf G}}(\omega)={\bf \Lambda}^{-\frac{1}{2}}\left[{\bf \Lambda}^{-\frac{1}{2}}\,{\bf K}{\bf \Lambda}^{-\frac{1}{2}}-{\bf 1}(\omega+i\epsilon)^2\right]^{-1}{\bf \Lambda}^{\frac{1}{2}}
\end{equation}
hence
\begin{equation}
\label{gftrue4}
{\hat {\bf G}}(\omega)={\bf \Lambda}^{-\frac{1}{2}}\,{\bf {\hat {\cal G}}}(\omega)\,{\bf \Lambda}^{\frac{1}{2}}
\end{equation}

${\hat {\bf G}}(\omega)$ is a non-symmetric matrix having the components
\begin{equation}
\label{gftrue5}
{\hat G}(\omega)_{pq}(\omega)= \xi^{-(p-q)}{\hat {\cal G}}_{pq}(\omega)
\end{equation}
where no sum over $p,q$ is performed. The same relation connects also the Green's functions of the time domain.
Let us evaluate now (\ref{greendf2}) in explicit form. To this end we represent
this Green's function in its spectral representation

\begin{equation}
\label{spectral1}
{\bf {\hat {\cal G}}}(\omega)=\sum_{m=0}^{N-1}\frac{1}{\omega_m^2-(\omega+i\epsilon)^2}{\bf v}^{(m)}\otimes{\bf v}^{(m)*}
\end{equation}
where $^{*}$ means complex conjugation and ${\bf v}^{(m)}$ denote the Bloch-eigenvectors of
${\bf L}$ and $({\bf a}\otimes {\bf b})_{pq}=a_pb_q$ indicates dyadic multiplication.
Hence
\begin{equation}
\label{spectral2}
{\hat {\cal G}}_{pq}(\omega)=\frac{1}{N}\sum_{m=0}^{N-1}\frac{e^{ik_m(p-q)}}{\omega_m^2-(\omega+i\epsilon)^2}
\end{equation}
Within the frequency band the Green's function (\ref{spectral2}) appears to be {\it complex} valued  due to the presence of the imaginary infinitesimal quantity $i\epsilon$. This infinitesimal term becomes irrelevant outside the frequency band where the Green's function is purely {\it real}.
In view of the symmetry of ${\hat {\cal G}}_{pq}$
one can conclude that ${\hat {\cal G}}_{pq}={\hat {\cal G}}_{|p-q|}$. To evaluate (\ref{spectral2}) explicitly we use the relation
\begin{equation}
\label{asrel1}
\frac{1}{N}\sum_{m=0}^{N-1} f(e^{ik_m})= \frac{1}{2\pi}\int_0^{2\pi}f(e^{ik})\,{\rm d}k
\end{equation}
which holds asymptotically for $N>>1$ where $k_m=\frac{2\pi m}{N}\in [0,2\pi]$ ($m=0,..,N-1$).
Hence by using (\ref{asrel1}) we can write the Green's function (\ref{spectral2}) with the explicit dispersion relation $\omega_m=\omega(k_m)$ of (\ref{eigenvals}) as the integral

\begin{equation}
\label{spectral3}
{\hat {\cal G}}_{pq}(\omega)=\frac{1}{2\pi}\int_0^{2\pi}\frac{e^{ik(p-q)}}{\omega^2(k)-(\omega+i\epsilon)^2}{\rm d}k
\end{equation}

To evaluate (\ref{spectral3}) explicitly by employing dispersion relation (\ref{eigenvals}) we set

\begin{equation}
\label{conv1}
\omega^2(k)-(\omega+i\epsilon)^2=
-\frac{\omega_0^2}{\xi}\mu^{-1}(\mu-e^{\phi})(\mu-e^{-\phi})
\end{equation}
with $\mu=e^{ik}$ and

\begin{equation}
\label{cosheq1}
a=\cosh{\phi}=\frac{(\Omega_D^2+\Omega_{0}^2-2\omega^2)}{(\Omega_D^2-\Omega_{0}^2)}
\end{equation}

where we introduced the lowest eigenfrequency\footnote{We emphasize that $\Omega_0=\frac{\omega_0}{\xi}|1-\xi|>0$ for $\xi\neq 1$.} $\Omega_0=\frac{\omega_0}{\xi}|1-\xi|$ and the highest (the Debye-) eigenfrequency $\frac{\omega_0}{\xi}(1+\xi)=\Omega_D$ of the chain.

The argument $\phi=\chi+i\psi$ in this equation is generally a {\it complex} quantity.
In the case that the zeros $e^{\pm\phi}$ of (\ref{conv1}) are located on the unit circle,
we have to replace in (\ref{cosheq1}) $\omega$ by $\omega+i\epsilon$. This is the case when $\omega$ is within the frequency band. To further evaluate (\ref{spectral3})
we write it as a complex integral over the unit circle
in the complex $\mu$-plane (where $\mu=e^{ik}$, ${\rm d}k=(i\mu)^{-1}{\rm d}\mu$)

\begin{equation}
\label{complexpint1}
{\hat {\cal G}}_{pq}=\frac{i\xi}{2\pi\omega_0^2}\oint_{|\mu|=1}\frac{\mu^{|p-q|}}{(\mu-e^{\phi})(\mu-e^{-\phi})}
{\rm d}\mu
\end{equation}

To evaluate this integral by using the theory of complex functions we have to determine the
residuum of the integrand within the unit circle $|\mu|=1$. Without loss of generality there is always one argument $\phi=\chi+i\psi$ with $ \chi\geq 0$ such that the pole at $e^{-\phi}$ is located {\it within the unit circle} since $|e^{-\phi}|=e^{-\chi}<1$. Hence evaluation of (\ref{complexpint1}) yields ${\hat {\cal G}}_{pq}=2\pi i Res(..)|_{\mu=e^{-\phi}}$ to arrive at

\begin{equation}
\label{complexpint1res}
{\hat {\cal G}}_{pq}=\frac{\xi}{2\omega_0^2}\frac{e^{-|p-q|\phi}}{\sinh{\phi}}
\end{equation}

To write the Green's function explicitly the following three cases have to be considered,
namely
\newline\newline\newline
\noindent {\bf Case (i)}: $\Omega_0=\frac{\omega_0}{\xi}|1-\xi|\leq \omega \leq \frac{\omega_0}{\xi}(1+\xi)=\Omega_D$
 \newline\newline
\noindent {\bf Case (ii)}: $0\leq \omega < \Omega_0$ \newline\newline
\noindent {\bf Case (iii)}: $\omega>\Omega_D$
\newline\newline\newline
\noindent {\bf Case (i)}: This case is relevant for the determination of the density of vibrational modes, i.e. when $\Omega_0\leq\omega\leq\Omega_D$ is within the band of eigenfrequencies.
In this case we observe that $a$ of eq. (\ref{cosheq1}) takes the values $a(\Omega_D)=-1\leq a \leq a(\Omega_0)= 1$. Hence we can rewrite (\ref{cosheq1}) as
\begin{equation}
\label{arel1cos}
a=\frac{(\Omega_D^2+\Omega_{0}^2-2\omega^2)}{(\Omega_D^2-\Omega_{0}^2)}=\cos{\varphi}\,,\hspace{2cm} 0\leq\varphi\leq \pi
\end{equation}
since $\varphi(\Omega_0)=0$ and $\varphi(\Omega_D)=\pi$.
In order to avoid the singularities of the integral (\ref{complexpint1}) being located on the complex unit circle we have to replace in (\ref{arel1cos}) $\omega$ by $\omega+i\epsilon$
with an infinitesimally small positive $\epsilon$
which is equivalent to replace $\varphi$ by $\varphi+i\epsilon$ in the limit $\epsilon$ tending to zero. Integral (\ref{complexpint1}) is then well defined, since the singularities
at $e^{\pm i\varphi}$ are shifted infinitesimally onto $e^{\pm i(\varphi+i\epsilon)}$.
By this regularization we have to put in (\ref{complexpint1res}) $\phi=\epsilon-i\varphi$ ($\chi=\epsilon\rightarrow 0+$) to arrive at

\begin{equation}
\label{complexpt3}
{\hat {\cal G}}_{pq}^{(i)}=\frac{i\xi}{2\omega_0^2}\frac{e^{i\varphi|p-q|}}{\sin{\varphi}}
\end{equation}
or
\begin{equation}
\label{complexpin4}
{\hat {\cal G}}_{pq}^{(i)}=\frac{i\xi}{2\omega_0^2}\frac{(a+i\sqrt{1-a^2})^{|p-q|}}{\sqrt{1-a^2}}
\end{equation}
where we have used that $\sin{\varphi}=+\sqrt{1-a^2}\geq 0$ in the interval $0\leq\varphi\leq \pi$ which assumes

\begin{equation}
\label{arel1sin}
\sin{\varphi}=\frac{2}{(\Omega_D^2-\Omega_{0}^2)}\sqrt{(\Omega_{D}^2-\omega^2)(\omega^2-\Omega_{0}^2)}
\end{equation}

With (\ref{arel1sin}) and (\ref{arel1cos}) we can rewrite the Green's function (\ref{complexpin4}) in the form
\begin{equation}
\label{complexfini}
{\hat {\cal G}}_{pq}^{(i)}(\omega)=\frac{i}{\sqrt{(\Omega_{D}^2-\omega^2)(\omega^2-\Omega_{0}^2)}}
\frac{\left(\Omega_{D}^2+\Omega_{0}^2-
2\omega^2+2i\sqrt{(\Omega_{D}^2-\omega^2)(\omega^2-\Omega_{0}^2)}\right)^{|p-q|}}{\left(\Omega_{D}^2-\Omega_{0}^2\right)^{|p-q|}}
\end{equation}
where (\ref{complexfini}) is restricted to case (i), i.e. to frequencies within the frequency band $\Omega_0\leq\omega\leq\Omega_D$. In this case we also obtain the density of vibrational modes by using the relation \cite{maradudin,mougios,3}
\begin{equation}
\label{densoscillators1}
\rho(\omega)=\frac{2\omega}{\pi}Im\,Tr\,{\hat {\cal G}}(\omega+i\epsilon)
\end{equation}

In (\ref{densoscillators1}) $Tr(..)$ denotes the trace of a matrix. We obtain for (\ref{densoscillators1})

\begin{equation}
\label{densoscil2}
\rho(\omega)=\frac{N\omega\xi}{\pi\omega_0^2\sin{\varphi}}=\frac{N\omega\xi}{\pi\omega_0^2\sqrt{1-a^2}}
\end{equation}
to arrive finally at
\begin{equation}
\label{oscden}
\rho(\omega)=\frac{2N\omega}{\pi \sqrt{(\Omega_{D}^2-\omega^2)(\omega^2-\Omega_{0}^2)}}\,,\hspace{1cm}
\Omega_{0} \leq \omega \leq \Omega_{D}
\end{equation}
which is defined and nonzero only within the frequency band and zero elsewhere
since the Green's function outside the band is purely real (cases (ii), (iii) below).

An important observation is that (\ref{oscden}) possesses two singularities, one at the Debye frequency $\Omega_D$ like in the case homogeneous linear chain, and
a second one at the lowest frequency $\Omega_0$ which disappears in the degenerate case
$\xi=1$ of a homogeneous chain. This second singularity at $\Omega_0$ is due to the suppression of all modes with frequencies $\omega <\Omega_0$ which are shifted to reappear at frequencies $\omega >\Omega_0$.

The mode density (\ref{oscden}) is normalized such that \cite{maradudin,3}

\begin{equation}
\label{norma}
\int_{\Omega_0}^{\Omega_D}\rho(\omega){\rm d}\omega=N
\end{equation}

To check this relation for (\ref{oscden}) we have to confirm that
\begin{equation}
\label{norma1}
\int_{\Omega_0}^{\Omega_D} \frac{2N\omega}{\pi \sqrt{(\Omega_{D}^2-\omega^2)(\omega^2-\Omega_{0}^2)}}  {\rm d}\omega=N
\end{equation}

It is convenient to use relation (\ref{arel1sin}) to evaluate

\begin{equation}
\label{norma2}
\int_{\Omega_0}^{\Omega_D}\rho(\omega)\,{\rm d}\omega=\frac{N\xi}{\pi\omega_0^2}\int_{\Omega_0}^{\Omega_D}
\frac{\omega}{\sin{\varphi(\omega)}}\,{\rm d}\omega
\end{equation}

To evaluate this integral we can directly make use of the substitution (\ref{arel1sin}) together with (\ref{arel1cos}) indicating that $\frac{4\omega}{\Omega_D^2-\Omega_0^2}{\rm d}\omega=\sin{\varphi}{\rm d}\varphi$.
Then we confirm the required relation (\ref{norma}), namely
\begin{equation}
\label{norma3}
\int_{\Omega_0}^{\Omega_D}\rho(\omega)\,{\rm d}\omega
=\frac{N}{\pi}(\varphi(\Omega_D)-\varphi{(\Omega_0)})=N
\end{equation}
where from (\ref{arel1cos}) follows that $\varphi(\Omega_0)=\arccos{(1)}=0$ and $\varphi(\Omega_D)=\arccos{(-1)}=\pi$.
The mode density (\ref{oscden}) degenerates in the case of a {\it homogeneous
linear chain represented by $\xi=1$} where $\Omega_{0}=0$ and $\Omega_{D}=\omega_D=2\omega_0$
to

\begin{equation}
\label{oschomchain}
\rho_{\xi=1}(\omega)=\frac{2N}{\pi\sqrt{(\omega_{D}^2-\omega^2)}}\,,\hspace{1cm}
0 \leq \omega \leq \omega_{D}
\end{equation}
which is the well known classical expression for the oscillator mode density of the homogeneous linear chain which can be found in any textbook of lattice dynamics (e.g. \cite{maradudin}). Expression (\ref{oschomchain}) can be also reobtained directly from the dispersion relation of the linear chain using $\rho_{\xi=1}(\omega)=2\frac{N}{2\pi}|\frac{dk}{d\omega}|$ and (\ref{eigenvals})\footnote{where the prefactor "2" accounts for the number of branches of the dispersion relation $\omega(k)$ in the $k$-space.} for $\xi=1$.
\newline\newline\newline
\noindent{\bf Case (ii)}:

In the range of frequencies $0\leq \omega < \Omega_0$, i.e. outside the band of eigenfrequencies the coefficient (\ref{cosheq1}) takes values $a>1$ thus we can put
\begin{equation}
\label{arel1cosh}
a=\frac{(\Omega_D^2+\Omega_{0}^2-2\omega^2)}{(\Omega_D^2-\Omega_{0}^2)}=\cosh{\chi}\,,\hspace{2cm} \chi>0
\end{equation}
where $\phi=\chi>0$ is purely real. Hence the Green's function (\ref{complexpint1res}) assumes in this case

\begin{equation}
\label{complexpint1res2}
{\hat {\cal G}}_{pq}^{(ii)}=\frac{\xi}{2\omega_0^2}\frac{e^{-|p-q|\chi}}{\sinh{\chi}}
\end{equation}
which is purely real valued.
In this expression $\sinh{\chi}=\sqrt{a^2-1}$ is the positive root because of $\chi>0$. We thus have

\begin{equation}
\label{arel1sinh}
\sinh{\chi}=\frac{2}{(\Omega_D^2-\Omega_{0}^2)}\sqrt{(\Omega_{D}^2-\omega^2)(\Omega_{0}^2-\omega^2)}
\end{equation}

Expression (\ref{complexpint1res2}) hence has to be rewritten as
\begin{equation}
\label{real2fini}
{\hat {\cal G}}_{pq}^{(ii)}(\omega)=\frac{1}{\sqrt{(\Omega_{D}^2-\omega^2)(\Omega_{0}^2-\omega^2)}}
\frac{\left(\Omega_{D}^2+\Omega_{0}^2-
2\omega^2-2\sqrt{(\Omega_{D}^2-\omega^2)(\Omega_{0}^2-\omega^2)}\right)^{|p-q|}}{\left(\Omega_{D}^2-\Omega_{0}^2\right)^{|p-q|}}
\end{equation}
where $0\leq\omega<\Omega_0$. We note that (\ref{real2fini}) contains the {\it static} limit ${\hat {\cal G}}_{pq}^{(ii)}(\omega=0)=({\bf L}^{-1})_{pq}$
which is obtained by putting $\omega=0$.
\newline\newline
\noindent{\bf Case (iii)}: It remains to discuss the case when the frequency exceeds the Debye frequency $\omega>\Omega_D$. In this case we have
\begin{equation}
\label{arel1cosh2}
a=\frac{(\Omega_D^2+\Omega_{0}^2-2\omega^2)}{(\Omega_D^2-\Omega_{0}^2)}=\cosh{(\chi+i\pi)}< -1\,,\hspace{2cm} \chi>0
\end{equation}
that is $a=-\cosh{\chi}$. Hence we have to put $\phi=\chi+i\pi$ in (\ref{complexpint1res}) to obtain the Green's function as
\begin{equation}
\label{rela3pint1res}
{\hat {\cal G}}_{pq}^{(iii)}=\frac{(-1)^{p-q+1}\xi}{2\omega_0^2}\frac{e^{-|p-q|\chi}}{\sinh{\chi}}
\end{equation}
Since $\chi>0$ the $\sinh{\chi}$ has to be as in case (ii) the positive root, namely
$\sinh{\chi}=\sqrt{a^2-1}$. Thus we have as in case (ii)

\begin{equation}
\label{arel3sinh}
\sinh{\chi}=\frac{2}{(\Omega_D^2-\Omega_{0}^2)}\sqrt{(\Omega_{D}^2-\omega^2)(\Omega_{0}^2-\omega^2)}
\end{equation}

To rewrite (\ref{rela3pint1res}) we have to note that $\cosh{\chi}=-a$ where $a$ is defined
by (\ref{arel1cosh2}). Taking this into account the Green's function has to be rewritten as
\begin{equation}
\label{real3fini}
{\hat {\cal G}}_{pq}^{(iii)}(\omega)=\frac{(-1)}{\sqrt{(\Omega_{D}^2-\omega^2)(\Omega_{0}^2-\omega^2)}}
\frac{\left(\Omega_{D}^2+\Omega_{0}^2-
2\omega^2+2\sqrt{(\Omega_{D}^2-\omega^2)(\Omega_{0}^2-\omega^2)}\right)^{|p-q|}}{\left(\Omega_{D}^2-\Omega_{0}^2\right)^{|p-q|}}
\end{equation}
where $\omega > \Omega_D$.

\section{The Continuum Limit}

We consider in this section the transition of system (\ref{hamiltonnth})
to a continuum. To this end we assume all particles having the same distance $h$ to each other with $h\rightarrow 0$ . We assume the chain has length $L=Nh$ where $L$ is kept finite and constant in this limiting process. As a consequence the total number of particles becomes infinite $N=\frac{L}{h}\rightarrow\infty$. Particle $p$ has then spatial coordinate $x=ph$ where
$p=0,..N-1$ varies $x$ quasi-continuous. The displacement $u_p(t)$ of particle $p$
has then to be replaced by the displacement field $u(x,t)$ where the the displacement of the
"neighbor" particle $p+1$ is given by $u(x+h,t)$.
In the continuum limit
any sum over $p$ can be replaced by an integral, namely
\begin{equation}
\label{sumrule}
\sum_{p=0}^{N-1}F_p\approx \int_0^L F(x)\frac{{\rm d}x}{h}
\end{equation}
where $L=Nh$ and ${\rm d}x\approx h$. This relation becomes asymptotically exact for a small $h\rightarrow 0$. We further introduce \cite{maugin}
\begin{equation}
\label{taylor}
u_{p+1}(t)-u_p(t)=u(x+h)-u(x)\approx h\frac{\partial u}{\partial x}(x,t)
\end{equation}
where we take into account only the linear order with respect to $h$.
The Hamiltonian (\ref{hamiltonnth})
can then be rewritten in this limiting case as the integral

\begin{equation}
\label{hamilcont}
H=\frac{m_0}{2h}\int_0^L e^{2\beta x}\left\{\left(\frac{\partial u}{\partial t}\right)^2
+h^2\omega_0^2\left(\frac{\partial u}{\partial x}\right)^2\right\}{\rm d}x
\end{equation}
where we have set $\xi^{2p}=e^{2\beta x}$ with $p=\frac{x}{h}$ and $\beta=\lim_{h\rightarrow 0}\frac{\ln{\xi}}{h}$. This is justified because we have to demand continuity of the function $\xi^p\approx \xi^{p+1}$, i.e. $\xi\approx 1$ when $h\rightarrow 0$. Hence it is justified to set in this limiting case
$\xi=e^{\beta h}\approx 1+\beta h$ where $\beta$ is constant.
In this way we can obtain the limit of $\xi^p=$ with $p=\frac{x}{h}$ for a given $x$
indeed in the form
\begin{equation}
\label{xip}
\xi^p=\lim_{h\rightarrow 0}(1+\beta h)^{\displaystyle \frac{x}{h}}=e^{\beta x}
\end{equation}

Since the integrand of (\ref{hamilcont}) represents a finite energy density
we have to demand in the limit $h\rightarrow 0$ the relations

\begin{equation}
\label{scaling}
m_0(h)=\rho_0h \,,\hspace{1.5cm}  \omega_0^2(h)=\frac{\Omega^2}{h^2} \,,\hspace{1.5cm} \beta=\frac{\ln{\xi}}{h}
\end{equation}
where we introduced the constants $\Omega$, $\rho_0$ and $\beta$. Then the energy density (energy per unit length) in (\ref{hamilcont}) can be written in the form
\begin{equation}
\label{Hamden}
{\cal H}= \frac{\rho_M(x)}{2}\left(\frac{\partial u}{\partial t}\right)^2+\frac{\mu(x)}{2}\left(\frac{\partial u}{\partial x}\right)^2
\end{equation}

The integrand represents the Hamiltonian density where the mass density $\rho_M(x)$ and elastic modulus $\mu(x)=\Omega^2\rho_M(x)$ are exponentially graded

\begin{equation}
\label{gradmod}
\rho_M(x)=\rho_0 \, e^{2\beta x}\,, \hspace{2cm} \mu(x)=\Omega^2 \rho_0\,e^{2\beta x}
\end{equation}

By conceiving integral (\ref{hamilcont}) as the energy functional, we obtain the equation of motion for
the unknown field $u(x,t)$ by employing the Hamiltonian principle \cite{goldstein} for fields, i.e. by putting the functional derivative $\frac{\delta H}{\delta u}=0$ to arrive at
the partial differential equation
\begin{equation}
\label{hamil2}
e^{2\beta x}\frac{\partial^2 u}{\partial t^2}=\Omega^2\frac{\partial}{\partial x}\left(e^{2\beta x}\frac{\partial u}{\partial x}\right)
\end{equation}
which corresponds to the case of absence of external forces.
This equation can be also obtained directly from (\ref{eq1}) by performing the transition by using (\ref{taylor}). Hence (\ref{hamil2}) is the continuous counterpart of (\ref{eq1}). In order to resolve (\ref{hamil2}) we make the transformation

\begin{equation}
\label{trafofield}
u(x,t)=y(x,t)e^{-\beta x}
\end{equation}
which is the continuous version of (\ref{subst1}). By this substitution (\ref{hamil2})
takes the form

\begin{equation}
\label{hamil4}
\frac{\partial^2 y}{\partial t^2}=
-{\cal L} y(x,t)
\end{equation}
where we introduced the differential operator

\begin{equation}
\label{hamil3diffop}
{\cal L}=-\Omega^2\left(\frac{\partial^2}{\partial x^2}-\beta^2\right)
\end{equation}

Equations (\ref{hamil4}) together with (\ref{hamil3diffop})
can be written in the form
\begin{equation}
\label{kleingordon}
\left[\frac{1}{\Omega^2}\frac{\partial^2 }{\partial t^2}-\frac{\partial^2 }{\partial x^2}\right]y+\beta^2y=0
\end{equation}
which has the form of a {\it Klein-Gordon} equation of the scalar field of a free moving (spin-0) particle \cite{davidov}. Further we introduce the periodic boundary conditions
\begin{equation}
\label{bc2cont}
y(x,t)=y(x+L,t)
\end{equation}
which are the continuous counterpart to (\ref{peroidboundc}).

We show now that (\ref{hamil4}) is the continuous counterpart of (\ref{newton2})
where the differential operator \newline ${\cal L}$ is nothing but the continuous counterpart of the matrix ${\bf L}$ of (\ref{matrixL1}). To this end we consider
\begin{equation}
\label{cont}
\begin{array}{lll}
\sum_q^{N-1}L_{pq}y_q(t)=\omega_0^2\left\{(1+e^{-2h\beta})y_p(t)-e^{-\beta h}(y_{p-1}(t)+y_{p+1}(t))\right\} && \nonumber \\
&& \nonumber \\
\hspace{2.6cm}=\omega_0^2\left\{(1-e^{-h\beta})^2y_p(t)-e^{-\beta h}(y_{p-1}(t)+y_{p+1}(t)-2y_p)\right\} &&
\end{array}
\end{equation}
Taking into account only terms up to order $h^2$ the first term in the last equation
contributes with $\beta^2h^2\omega_0^2y(x,t)=\Omega^2\beta^2y(x,t)$ and the second term with
$-h^2\omega_0^2\frac{\partial^2}{\partial x^2}y(x,t)=-\Omega^2\frac{\partial^2}{\partial x^2}y(x,t)$ to arrive finally at

\begin{equation}
\label{limform2}
\sum_q^{N-1}L_{pq}y_q(t)\approx -\Omega^2(\frac{\partial}{\partial x^2}-\beta^2)y(x,t)={\cal L}y(x,t)
\end{equation}

The normalized eigenfunctions of operator ${\cal L}$ can be written as the continuous
version of the Bloch eigenvectors (\ref{blocheig})
\begin{equation}
\label{blocheigcont}
v^{(m)}(x)=\frac{1}{\sqrt{L}} e^{iK_m x}\,,\hspace{2cm} K_m=\frac{2\pi}{L}m
\end{equation}
with $m=0, \pm 1, \pm 2, ..,\pm \infty$. The eigenfunctions $v^{(m)}(x)$ fulfil the periodic boundary conditions (\ref{bc2cont}). The dispersion relation is then obtained by the eigenvalues of ${\cal L}$, namely

\begin{equation}
\label{eigval}
{\cal L}v^{(m)}(x)=\omega^2(K_m)v^{(m)}(x)
\end{equation}
and yields
\begin{equation}
\label{disp}
\omega^2(K_m)=\Omega^2(K_m^2+\beta^2)
\end{equation}

This dispersion relation has to coincide with the dispersion relation
(\ref{eigenvals}) in the continuum limit. To this end we rewrite
(\ref{eigenvals}) where we note the relation between the wave numbers $k_m=hK_m$ and obtain

\begin{equation}
\label{eigenvalsrew}
\omega_m^2=\omega_0^2(e^{iK_mh-\beta h}-1)(e^{-iK_mh-\beta h}-1)\approx \omega_0^2h^2(K_m^2+\beta^2)=\Omega^2(K_m^2+\beta^2)
\end{equation}
which recovers (\ref{disp}). We note that the dispersion relation (\ref{eigenvalsrew})
represents a discrete set of eigenfrequencies due to the discrete values of $K_m$. This is a consequence of the periodic boundary conditions (\ref{bc2cont}) with a finite length
$L$ of the chain.
In the limiting case $L\rightarrow\infty$ the spectrum of eigenfrequencies becomes quasi-continuous.
\newline\newline\newline
\noindent {\bf Continuum limit of the Green's function}. In the following we deduce the continuous counterpart of Green's function (\ref{complexpint1res}) and above cases (i)-(iii).
Then we have $\Omega_0=\frac{|1-\xi|}{\xi}\omega_0=\omega_0(e^{\beta h}-1)e^{-\beta h}\approx\beta \omega_0 h=\beta\Omega$ and in the same way $\Omega_D\approx 2\frac{\Omega}{h}+\beta\Omega\rightarrow\infty$. This observation, that the lowest frequency
is at $\Omega_0=\beta\Omega$ with no finite Debye frequency is also reflected by the dispersion relation (\ref{eigenvalsrew}). We have hence in the continuum limit only the two cases:
\newline\newline\newline
\noindent {\bf Case (i)}: $\Omega_0=\beta\Omega\leq \omega < \infty $ \newline\newline
\noindent {\bf Case (ii)}: $0\leq \omega < \Omega_0=\beta\Omega$
\newline\newline
To consider these two cases it is useful to write coefficient (\ref{cosheq1}) in this approximation
\begin{equation}
\label{coefapp}
a=1+\frac{h^2}{2}(\beta^2-\frac{\omega^2}{\Omega^2})
\end{equation}
where the lowest non-vanishing order in $h$ is the second order term which has been taken into account. In view of (\ref{coefapp}) it follows that in case (i) $a<1$ and in case (ii) $a>1$.
\newline\newline
\noindent {\bf Case (ii)}: Here we can put
\begin{equation}
\label{coshlim}
a=\cosh{\chi}=1+\frac{h^2}{2}(\beta^2-\frac{\omega^2}{\Omega^2})>1\,\hspace{2cm} \chi > 0
\end{equation}
since $\beta^2-\frac{\omega^2}{\Omega^2}\geq 0$ and by taking into account only the dominating order in $h$ we have
\begin{equation}
\label{sinhapprox}
\sinh{\chi}=h\sqrt{\beta^2-\frac{\omega^2}{\Omega^2}}
\end{equation}
where the positive root has to be taken since $\chi>0$. Let us consider the expression (\ref{complexpint1res2}) for ${\hat {\cal G}}_{pq}^{(ii)}$ in the continuous approximation. Here we have to put $p-q=\frac{x_1-x_2}{h}=\frac{x}{h}$ and the fact that
a matrix multiplication with the Green's function becomes a convolution integral where
the matrix of the Green's function becomes a convolution kernel.
Since the summation of a matrix multiplication has to be transfered into an integral
by taking into account (\ref{sumrule}), we have to multiply the discrete Green's function with
a factor $\frac{1}{h}$ to arrive at ${\hat g}(|x|,\omega)=\lim_{h\rightarrow 0}\left(\frac{1}{h}{\hat {\cal G}}_{pq}^{(ii)}\right)$. Then we can generate the continuous limit of the Green's function by constructing $e^{-\chi}$
from (\ref{coshlim}) and (\ref{sinhapprox}) and by using
the Green's function matrix of (\ref{complexpint1res2}) in the form

\begin{equation}
\label{complexlim}
{\hat g}^{(ii)}(|x|,\omega)=\lim_{h\rightarrow 0} \frac{1}{h}\frac{h^2e^{-\beta h}}{2\Omega^2}\frac{\left(1+\frac{h^2}{2}(\beta^2-\frac{\omega^2}{\Omega^2})-
h\sqrt{\beta^2-\frac{\omega^2}{\Omega^2}}\right)^{\frac{|x|}{h}}}
{h\sqrt{\beta^2-\frac{\omega^2}{\Omega^2}}}
\end{equation}
which simplifies when we take into account the relation
\begin{equation}
\label{hrela}
\lim_{h\rightarrow 0}\left(1+c_1 h +c_2h^2 \right)^{\frac{x}{h}}=e^{c_1x}
\end{equation}
where the quadratic term in $c_2h^2$ and higher powers in $h$ are irrelevant. Utilizing this relation yields the Green's function in the form

\begin{equation}
\label{greenffiilim}
{\hat g}^{(ii)}(|x|,\omega)=\frac{1}{2\Omega^2}\frac{e^{-|x|
\sqrt{\beta^2-\frac{\omega^2}{\Omega^2}}}}
{\sqrt{\beta^2-\frac{\omega^2}{\Omega^2}}}
\end{equation}
where $0\leq \omega < \beta\Omega$. It is important to note that because of the periodic boundary conditions (\ref{bc2cont}) which are to be fulfilled, the Green's function has the form (\ref{greenffiilim}) in the interval $0\leq |x| < L$ and has to be $L$-periodically
continued. In the limiting case $L\rightarrow\infty$ (\ref{greenffiilim}) is defined
for all $0\leq |x|< \infty$.
\newline\newline
{\bf Case (i)}: In a similar way we generate the continuous limit in the band of eigenfrequencies, i.e. when $\beta\Omega \leq \omega<\infty$.
In this case we have $a<1$ and put
\begin{equation}
\label{coefappcos}
a=\cos{\varphi}=1-\frac{h^2}{2}(\frac{\omega^2}{\Omega^2}-\beta^2)<1\,,\hspace{2cm}\varphi > 0
\end{equation}
and
\begin{equation}
\label{sinapprox}
\sin{\varphi}=h\sqrt{\frac{\omega^2}{\Omega^2}-\beta^2}
\end{equation}
where the positive root has to be taken since $\varphi > 0$.
In a complete analogous way as above in (\ref{complexlim}) we consider the limit $h\rightarrow 0$ of (\ref{complexpt3}) by putting ${\hat g}=\frac{1}{h}{\hat {\cal G}}_{pq}^{(i)}$
and obtain

\begin{equation}
\label{complexpt3contlim}
{\hat g}^{(i)}(|x|,\omega)=\frac{i}{2\Omega^2}\frac{e^{i|x|
\sqrt{\frac{\omega^2}{\Omega^2}-\beta^2}}}
{\sqrt{\frac{\omega^2}{\Omega^2}-\beta^2}}
\end{equation}
where $\omega\geq\beta\Omega$. By utilizing (\ref{complexpt3contlim}) we can also generate
the continuous limit of the mode density which is non-zero within the band of eigenfrequencies $\omega\geq\beta\Omega$. To this end we have to employ the continuous counterpart of
(\ref{densoscillators1}) to arrive at

\begin{equation}
\label{densoscillators1cont}
\rho(\omega)=\frac{2\omega L}{\pi}Im\,{\hat g}^{(i)}(|x|=0,\omega)
\end{equation}
which yields

\begin{equation}
\label{denmode}
\rho(\omega)=\frac{L\omega}{\pi\Omega^2\sqrt{\frac{\omega^2}{\Omega^2}-\beta^2}}
\end{equation}

As in the discrete case this mode density is only defined within the band $\omega\geq\beta\Omega$ and zero outside for $0\leq\omega<\beta\Omega$. Expression (\ref{denmode}) is also in accordance with the continuous limit of (\ref{oscden}), by using the asymptotic relations $\Omega_0=\beta\Omega$ and $\Omega_D= \frac{2\Omega}{h}$.

The deduced Green's function ${\hat g}(|x|,\omega)$ of equations (\ref{greenffiilim})
and (\ref{complexpt3contlim}) is the fundamental solution of the frequency domain
representation of (\ref{hamil4}) and fulfills the equation

\begin{equation}
\label{greendeflim}
({\cal L}-\omega^2) g(|x|,\omega)=\delta(x)
\end{equation}
where $\delta(x)$ denotes the Dirac $\delta$-function. This equation is noting but the continuous analogue of the matrix equation leading to (\ref{greendf2}) in the discrete case.
We note that when we impose periodic boundary conditions the Green's function (\ref{greenffiilim})
and (\ref{complexpt3contlim}) is defined at $0\leq |x| <L$ and has to be $L$-periodically continued. In the case $L\rightarrow\infty$ of an infinitely long chain these expressions hold for the entire range $0\leq |x| <\infty$ and represent the Green's function of the 1D infinite space. Equation (\ref{greendeflim}) can be rewritten in the form of an inhomogeneous Helmholtz equation, namely
\begin{equation}
\label{hamil45}
-\Omega^2\left(\frac{d^2}{d x^2}+\frac{\omega^2}{\Omega^2}-\beta^2\right)g(|x|,\omega)=\delta(x)
\end{equation}
Obviously the character of solution of this equation depends on the sign
of $\frac{\omega^2}{\Omega^2}-\beta^2$ reflecting the two cases (i) and (ii).
\newline
\noindent It remains to verify that expressions (\ref{greenffiilim}) and (\ref{complexpt3contlim}) indeed solve (\ref{hamil45}).
To this end we consider
\begin{equation}
\label{veri1}
\frac{d^2}{dx^2}ce^{\lambda|x|}=\lambda^2 c\,e^{\lambda|x|}+2c\, \lambda \, \delta(x)
\end{equation}
where $\Theta(x)$ denotes the Heaviside unit step function and $c$ denotes a constant. We further have used
\begin{equation}
\label{montant}
\frac{d^2|x|}{dx^2}=\frac{d}{dx}(\Theta(x)-\Theta(-x))=2\delta(x)
\end{equation}
and $\delta(x)e^{\lambda|x|}=\delta(x)$.
We hence can rewrite (\ref{veri1}) in the form
\begin{equation}
\label{reveri1}
(\frac{d^2}{dx^2}-\lambda^2)ce^{\lambda|x|}=2c\lambda \delta(x)
\end{equation}

By putting $c=\frac{-1}{2\lambda\Omega^2}$ we recover with $\lambda=i\sqrt{\frac{\omega^2}{\Omega^2}-\beta^2}$ expression (\ref{complexpt3contlim}) of case (i) and
with $\lambda=-\sqrt{\beta^2-\frac{\omega^2}{\Omega^2}}$ expression (\ref{greenffiilim}) of case (ii). We note that solutions with $\lambda$'s of opposite signs also solve (\ref{hamil45}), however they correspond to physically inadmissible\footnote{$\lambda= -i\sqrt{\frac{\omega^2}{\Omega^2}-\beta^2}$ leads to a non-causal Green's function in the time domain in case (i) and $\lambda=\sqrt{\beta^2-\frac{\omega^2}{\Omega^2}}$ to a divergent solution for $L\rightarrow\infty$ at $x\rightarrow\infty$ in case (ii).} solutions of (\ref{hamil45}). The continuous counterpart to (\ref{gftrue5}), i.e. the Green's function of the "true" displacement field $u(x,t)$ is
given by
\begin{equation}
\label{gftrue5contlim}
{\tilde g}(|x|,\omega)= e^{-\beta x}{\hat g}(|x|,\omega)
\end{equation}
where ${\hat g}(|x|,\omega)$ is the Green's function given by (\ref{greenffiilim}) and (\ref{complexpt3contlim}), respectively. We further observe that (\ref{gftrue5contlim}) fulfills

\begin{equation}
\label{fonctiondegreeneq}
-\Omega^2\left\{\frac{d}{dx}\left(e^{2\beta x}\frac{d {\tilde g}}{dx}\right)+\frac{\omega^2}{\Omega^2}e^{2\beta x}{\tilde g}\right\}=\delta(x)
\end{equation}
which is the continuous counterpart of the matrix equation defining the discrete Green's function (\ref{gftrue}).

\section{Conclusions}

Expressions (\ref{complexfini}), (\ref{real2fini}), (\ref{real3fini}), respectively, provide the explicit
representations for the dynamic $N\times N$ Green's function defined by (\ref{greendf2})
for the entire range of frequencies. This Green's function is linked with the Green's function of the "true" displacements simply by (\ref{gftrue5}).
The degenerate case $\xi=1$ recovers the well known expressions for the Green's function of the homogeneous linear chain which have been known for a long time (see e.g. \cite{maradudin,3}).

We also deduced the continuous limit of this system which corresponds to a 1D graded linear elastic material. The transformed equations of motion take in the continuous limit the form
of a Klein-Gordon equation (\ref{kleingordon}) where the limiting transformation of the Green's function
yields the expressions (\ref{greenffiilim}) and (\ref{complexpt3contlim}), respectively,
which solve (\ref{hamil45}).

In spite of the simplicity of the model system, the exponentially graded linear chain and its continuous counterpart seem to be one of the simplest systems representing one-dimensionally graded materials which are accessible to a concise analysis.

The band of natural frequencies, i.e. its width as well as the highest and lowest frequencies can be controlled simply by the grading parameter $\xi$. This suggests the possibility to design exponentially graded FGM with frequency bands according to desired dynamic properties.

Possible further extensions of the present analysis could be on
linear chains with {\it self-similar} (fractal) properties
\cite{michelarxiv}. We hope that the present analysis inspires also
work in such directions.

\section{Acknowledgements}
We dedicate this paper to our dear friend Valery M. Levin on the occasion of his birthday
and wish him many years of creative power.

\end{document}